\documentclass{amsart}
\usepackage[final]{microtype}
\usepackage[a4paper,text={160mm,254mm},centering]{geometry}
\usepackage{hyperref,amsfonts,amssymb,amsthm,mathtools,xcolor,lmodern,mathrsfs}
\usepackage{tikz-cd}
\tikzcdset{arrow style=tikz}

\renewcommand{\setminus}{\mathbin{\rule[0.2em]{0.67em}{0.12em}}}%
\renewcommand{\mathbb}{\mathbf}

\renewcommand{\leq}{\leqslant}
\renewcommand{\geq}{\geqslant}
\renewcommand{\subset}{\subseteq}

\usepackage[doi=false,isbn=false,sorting=nyt,url=false,eprint=false,backref=true, style=numeric-comp]{biblatex}

\theoremstyle{definition}
\newtheorem{remark}{Remark}[section]
\newtheorem{definition}[remark]{Definition}
\newtheorem{construction}[remark]{Construction}

\theoremstyle{plain}
\newtheorem{theorem}[remark]{Theorem}

\newtheorem{lemma}[remark]{Lemma}
\newtheorem{proposition}[remark]{Proposition}

\numberwithin{equation}{remark}
\newtheoremstyle{void} {.5\baselineskip plus .2\baselineskip minus .2\baselineskip} {.5\baselineskip plus .2\baselineskip minus .2\baselineskip} {\normalfont} {} {\bfseries} {.}
{5pt plus 1pt minus 1pt} 
{\thmname{#1 }\thmnumber{#2}\thmnote{.~{#3}}}
\theoremstyle{void}

\numberwithin{equation}{remark}

\title{On an elementary proof of Weil II}
\date{November 22, 2026}
\author{Dingxin Zhang}
\begin{document}
\maketitle

\section{Introduction}

In \cite{deligne_weil2}, Deligne developed a fundamental theory of
weights that has since found numerous applications in number theory
and arithmetic geometry.  Recently, L.~Alpöge suggested a proof of
Deligne's theorem in footnote~3 of his
\href{https://alpo.ge/remark.pdf}{remark} on joint work with Claude on
the Riemann hypothesis.  An expanded account of this finite-field
argument appeared on
\href{https://hexagonmath.org/2611.00789v1}{\emph{Hexagon}}, uploaded
by F.~Alderliefste.  As that account is AI-generated and not
well-polished, we offer this note in the hope of making the argument
clearer.  \emph{We make no claim to originality.}

To avoid any misunderstanding, we should clarify that ``elementary''
here refers only to the analytic core of the argument.  The formalism
of \(\ell\)-adic cohomology, including the trace formula, Artin's
vanishing theorem, and Deligne's perverse weak Lefschetz theorem, is
still needed.

The main result is the following theorem.

\begin{theorem}
\label{theorem:weil2}
Assume \(\ell\nmid q\). Let \(\mathcal{E}\) be a lisse
\(\overline{\mathbb{Q}}_{\ell}\)-sheaf on
\(\mathbb{A}^1_{\mathbb{F}_q}\) of pointwise weight \(\leq w\) for an embedding
\(\iota\colon\overline{\mathbb{Q}}_{\ell}\hookrightarrow\mathbb{C}\).
Then every Frobenius eigenvalue \(\alpha\) on
\(\mathrm{H}_c^1(\mathbb{A}^1_{\overline{\mathbb{F}}_q};\mathcal{E})\) satisfies
\[
|\iota(\alpha)|\leq q^{\frac{w+1}{2}}.
\]
\end{theorem}

Here, saying that \(\mathcal{E}\) has pointwise weight \(\leq w\) means that
for any \(x \in \mathbb{F}_{q^{m}}\) and any eigenvalue \(\gamma\) of the
geometric Frobenius \(F_{x}\) at \(x\), we have
\begin{equation*}
|\iota(\gamma)| \leq q^{\frac{mw}{2}}.
\end{equation*}

\begin{remark}
Theorem~\ref{theorem:weil2} immediately implies the corresponding
statement with \(\mathbb{A}^1_{\mathbb{F}_q}\) replaced by any
nonempty open subscheme
\(U = \mathbb{A}^1_{\mathbb{F}_q} \setminus \{h = 0\}\), where
\(h \in \mathbb{F}_q[t]\) is nonconstant. Indeed, the ``Abhyankar
map'' (cf.~\cite[Lecture 1, ``First reduction\ldots'']{katz_weil-2}
and \cite[Theorem~1]{abhyankar_covering-of-algebraic-curves})
\(\pi\colon U \to \mathbb{A}^1_{\mathbb{F}_q}\) given by
\(t \mapsto t + h(t)^{-q}\) is finite étale, and we can apply
Theorem~\ref{theorem:weil2} to \(\pi_{\ast}\mathcal{E}\).

Alternatively, the proof in this note adapts directly to any nonempty
open \(U \subset \mathbb{A}^1_{\mathbb{F}_q}\) with minimal changes
(replacing \(\mathbb{A}^n\) by \(U^n\) and extending the trace
function by \(0\) outside \(U\)).  We restrict to \(\mathbb{A}^1\) for
brevity.
\end{remark}

\begin{remark}
By the preceding remark and Laumon~\cite[(4.3.2.1)]{laumon_fourier-transform},
Theorem~\ref{theorem:weil2} implies the purity of the intersection
complex in dimension \(1\)
\cite[Théorème~2]{deligne_weil2}. From this,
together with a more-or-less standard dévissage argument, Deligne
deduced his general theorem on the weights of extraordinary direct
images, see \cite[Théorème~3.3.1]{deligne_weil2}.

In fact, the \href{https://hexagonmath.org/2611.00789v1}{Hexagon note}
proved Deligne's general theorem directly, without first reducing to
the case of curves. However, we have chosen to treat only the
\(\mathbb{A}^1\)-case: it already contains the heart of the argument,
exhibits the key technical ideas much more clearly, and keeps the
notation lighter. The reader should easily be able to adapt the proof
to the general setting.  Alternatively, for the dévissage from the
\(\mathbb{A}^1\)-case, we refer to Deligne's article
\cite{deligne_weil2}, the textbook
\cite{kiehl-weissauer_weil-conjectures-perverse-sheaves}, or Danilov's
survey \cite[p.~116]{danilov_cohomology-of-varieties}.
\end{remark}

\begin{remark}
Theorem~\ref{theorem:weil2} has hitherto had three different proofs:
Deligne's original approach based on his interpretation of the method
of Hadamard and de la Vallée Poussin, Laumon's approach via the
Fourier transform \cite{laumon_fourier-transform}, and a variant of
Laumon's proof by Katz \cite{katz_weil-2} using Larsen's
alternative. All of these essentially rely on the so-called tensor
power trick.
The approach in this note is no different: we deduce the theorem from
the following a priori weaker result.
\end{remark}

\begin{proposition}\label{proposition:main}
For each \(n\geq 2\), the group
\(\mathrm{H}_c^n(\mathbb{A}^n;\mathcal{E}^{\boxtimes n})\)
has weights \(\leq n(w+1)+3\).
\end{proposition}

Here and throughout, \emph{cohomology is taken over
  \(\overline{\mathbb{F}}_q\)}, and saying that a cohomology group has
weights \(\leq w\) means that every eigenvalue \(\beta\) of the
geometric Frobenius \(F\) on it satisfies
\(|\iota(\beta)| \leq q^{\frac{w}{2}}\).

\begin{proof}[Proof that Proposition~\ref{proposition:main} \(\Rightarrow\) Theorem~\ref{theorem:weil2}]
By the Künneth formula, there is a Frobenius-equivariant injection
\[
\mathrm{H}_c^1(\mathbb{A}^1;\mathcal{E})^{\otimes n} \to
\mathrm{H}^n_c(\mathbb{A}^n ; \mathcal{E}^{\boxtimes n}).
\]
If \(\alpha\) is a Frobenius eigenvalue of
\(\mathrm{H}^1_c(\mathbb{A}^1; \mathcal{E})\), then \(\alpha^n\) is a
Frobenius eigenvalue of
\(\mathrm{H}^n_c(\mathbb{A}^n; \mathcal{E}^{\boxtimes n})\). The
proposition therefore gives
\[
|\iota(\alpha)|^n\leq q^{\frac{n(w+1)+3}{2}},
\quad\text{so}\quad
|\iota(\alpha)|\leq q^{\frac{w+1}{2}+\frac{3}{2n}}.
\]
Letting \(n\) tend to infinity proves the theorem.
\end{proof}

In the remainder of this note, we prove the proposition. Note that if
we pass to a finite extension of \(\mathbb{F}_{q}\) of degree \(d\),
both the size of the ground field and the absolute values of the
Frobenius eigenvalues are raised to the \(d\)-th power, so taking
\(d\)-th roots recovers the bound over \(\mathbb{F}_{q}\). In what
follows, we will therefore freely pass to finite extensions of
\(\mathbb{F}_{q}\).

\section{Proof of Proposition~\ref{proposition:main}}

We shall prove the proposition via Deligne's
\emph{perverse weak Lefschetz theorem} together with a power-saving
\emph{exponential sum estimate}, Lemma~\ref{lemma:key}.

We begin with some general remarks.
Suppose \(\mathcal{G}\) is a \emph{lisse}
\(\overline{\mathbb{Q}}_{\ell}\)-sheaf on \(\mathbb{A}^{n}\). Then
for any closed embedding of a smooth hypersurface
\begin{equation*}
i\colon Y \to \mathbb{A}^{n}
\end{equation*}
with complement
\(j\colon \mathbb{A}^{n} \setminus Y \to \mathbb{A}^{n}\), the
\emph{relative purity theorem}
\cite[Exposé XVI, 3.7, 3.8, and 3.10(a)]{sga4}
states that there is a natural isomorphism
\begin{equation*}
i^{!} \mathcal{G} \simeq i^{\ast} \mathcal{G}[-2](-1).
\end{equation*}
Therefore, applying \(R\Gamma_{c}(\mathbb{A}^{n},-)\) to the homotopy
fiber sequence
\begin{equation*}
i_{\ast}i^{!}\mathcal{G} \to \mathcal{G} \to Rj_{\ast}j^{\ast}\mathcal{G}
\end{equation*}
and taking cohomology, we obtain a long exact sequence
\begin{equation}\label{eq:les}\tag{†}
0 \to \mathrm{H}^{n}_{c}(\mathbb{A}^{n}; \mathcal{G})
\to \mathrm{H}^{n}_{c}(\mathbb{A}^{n}; Rj_{\ast}j^{\ast}\mathcal{G})
\to \mathrm{H}^{n-1}_{c}(Y; i^{\ast}\mathcal{G}(-1)) \xrightarrow{\mathrm{Gysin}} \mathrm{H}^{n+1}_{c}(\mathbb{A}^{n}; \mathcal{G}) \to \cdots.
\end{equation}
Here, the exact sequence begins at \(\mathrm{H}^{n}_{c}\) thanks to
Artin's vanishing theorem and Poincaré duality (\cite[Exposé~XIV, 3.1]{sga4} or \cite[Théorème 4.1.1]{bbd}).

In the situation of Proposition~\ref{proposition:main},
examining the exact sequence \eqref{eq:les} with
\(\mathcal{G} = \mathcal{E}^{\boxtimes n}\),
we see that the proposition follows at once
if we can prove that the weights of
\(\mathrm{H}^n_c(\mathbb{A}^n;Rj_{\ast}j^{\ast}\mathcal{E}^{\boxtimes n})\)
are \(\leq n(w+1)+3\) for \emph{one} smooth hypersurface \(Y\).
From the discussion above, we deduce the following criterion.

\begin{lemma}\label{lem:trace-bound}
If there exist a smooth hypersurface \(i\colon Y \to \mathbb{A}^{n}\)
and a constant \(C > 0\) such that
\begin{equation*}
\bigl|\iota\operatorname{Tr}\bigl(F^{m}, \mathrm{H}^{n}_{c}(\mathbb{A}^{n}; Rj_{\ast}j^{\ast}\mathcal{E}^{\boxtimes n})\bigr)\bigr| \leq C (q^{m})^{\frac{n(w+1)+3}{2}}
\end{equation*}
holds for all \(m \geq 1\), then Proposition~\ref{proposition:main} holds.
\end{lemma}

\begin{proof}
The transition from traces to eigenvalues is standard.  For the
convenience of the reader, we recapitulate the argument.  Let \(\beta_j\)
be the Frobenius eigenvalues, with multiplicity, on
\(\mathrm{H}_c^n(\mathbb{A}^n;Rj_{\ast}j^{\ast}\mathcal{E}^{\boxtimes
  n})\).  By hypothesis,
\[
\left|\sum_j\iota(\beta_j)^m\right|
\leq C q^{\frac{m(n(w+1)+3)}{2}}
\]
for all \(m\geq 1\). Consequently, the power series
\[
\sum_{m\geq 1}\left(\sum_j\iota(\beta_j)^m\right)T^{m-1}
=\sum_j\frac{\iota(\beta_j)}{1-\iota(\beta_j)T}
\]
converges for \(|T|<q^{-\frac{n(w+1)+3}{2}}\). Since each
\(\iota(\beta_j)^{-1}\) (for \(\beta_j\neq 0\)) is a pole of the
right-hand side, we obtain
\(|\iota(\beta_j)|\leq q^{\frac{n(w+1)+3}{2}}\) for all \(j\).
By the injection in \eqref{eq:les}, the same bound holds for
\(\mathrm{H}_c^n(\mathbb{A}^n;\mathcal{E}^{\boxtimes n})\), which proves
the proposition.
\end{proof}

It turns out that the cohomology spaces
\(\mathrm{H}_{c}^{\ast}(\mathbb{A}^{n};Rj_{\ast}j^{\ast}\mathcal{E}^{\boxtimes n})\)
are closely related to exponential sums.
To explain this, fix an integer \(m \geq 1\) and introduce some notation.
\begin{itemize}
\item For \(x\in\mathbb{F}_{q^m}\), let
\(t(x)=\iota\operatorname{Tr}(F_x, \mathcal{E}_{\overline{x}})\),
where \(F_x\) is geometric Frobenius at the \(\mathbb{F}_{q^m}\)-point \(x\).
\item Fix a nontrivial additive character
\(\psi\colon\mathbb{F}_q\to\mathbb{C}^\times\).
Let
\(\psi_m=\psi\circ\operatorname{Tr}_{\mathbb{F}_{q^m}/\mathbb{F}_q}\).
\end{itemize}
With this notation, we have the following lemma.

\begin{lemma}\label{lem:cayley-trace}
Suppose \(Y\) is a smooth hypersurface in \(\mathbb{A}^{n}\) defined
by \(P = 0\).  For geometric Frobenius \(F\) over \(\mathbb{F}_q\), we
have
\[
\sum_{e} (-1)^{e}\iota\operatorname{Tr}\bigl(F^m,
\mathrm{H}_c^e(\mathbb{A}^n;Rj_*j^*(\mathcal{E}^{\boxtimes n}))\bigr)
= -\sum_{h\in\mathbb{F}_{q^m}^{\times}}
\sum_{x\in\mathbb{F}_{q^m}^n}
t(x_1)\cdots t(x_n)\psi_m(hP(x)).
\]
\end{lemma}

\begin{proof}
By the Grothendieck trace formula for \(\mathbb{A}^{n}\) and \(Y\)
\cite[Rapport, Théorème 3.2]{sga4.5}, and the exact
sequence \eqref{eq:les}, we have
\begin{equation}\label{eq:difference}
\sum_{e}(-1)^e\iota\operatorname{Tr}\bigl(F^m,
\mathrm{H}_c^e(\mathbb{A}^n;Rj_*j^*(\mathcal{E}^{\boxtimes n}))\bigr)
=
\sum_{x\in\mathbb{F}_{q^m}^n}t(x_1)\cdots t(x_n)
-q^m\sum_{x\in Y(\mathbb{F}_{q^m})}t(x_1)\cdots t(x_n).
\end{equation}
The factor \(q^m\) arises from the Tate twist \((-1)\). By
orthogonality of characters,
\[
\sum_{h\in\mathbb{F}_{q^m}}\psi_m(hP(x))=
\begin{cases}
q^m,& P(x)=0,\\
0,& P(x)\ne0.
\end{cases}
\]
Hence
\[
q^m\sum_{x\in Y(\mathbb{F}_{q^m})}t(x_1)\cdots t(x_n)
=\sum_{h\in\mathbb{F}_{q^m}}\sum_{x\in\mathbb{F}_{q^m}^n}
t(x_1)\cdots t(x_n)\psi_m(hP(x)).
\]
The term \(h=0\) on the right equals the sum over \(\mathbb{A}^n\) in
the trace formula above. Substituting this into \eqref{eq:difference}
completes the proof.
\end{proof}

As we will see from Lemma~\ref{lemma:key} shortly, when \(P\) is
chosen suitably, the right-hand side of the identity in
Lemma~\ref{lem:cayley-trace} admits an elementary power-saving
estimate. The issue is the left-hand side. To apply
Lemma~\ref{lem:trace-bound}, we need to bound the Frobenius traces on
\(\mathrm{H}^{n}_{c}(\mathbb{A}^{n};Rj_{\ast}j^{\ast}\mathcal{E}^{\boxtimes
  n})\); however, the left-hand side is an alternating sum, and
potential cancellation may prevent us from isolating the degree-\(n\)
term. To this end, we wish to arrange for all terms with
\(e \neq n\) to vanish. This is where Deligne's perverse weak
Lefschetz theorem enters.

\begin{definition}
Let \(Y \subset \mathbb{A}^{n}\) be a hypersurface.  We say that the
\emph{weak Lefschetz theorem} holds for the pair \((Y,\mathcal{G})\)
if
\(\mathrm{H}^{e}_{c}(\mathbb{A}^{n}; Rj_{\ast}j^{\ast}\mathcal{G}) =
0\) for all \(e \neq n\).

When \(Y\) is smooth, weak Lefschetz holds for \((Y,\mathcal{G})\) if
and only if the Gysin map
\begin{equation*}
\mathrm{H}_{c}^{n-1+e}(Y; i^{\ast}\mathcal{G}(-1)) \to \mathrm{H}^{n+1+e}_{c}(\mathbb{A}^{n};\mathcal{G})
\end{equation*}
is surjective for \(e = 0\) and bijective for \(e > 0\).
\end{definition}

Deligne proved a perverse weak Lefschetz theorem in \cite[Corollary
A.5]{katz_affine-cohomological-transforms}.  The special case we need
is formulated below.

\begin{theorem}
\label{theorem:deligne-weak-lefschetz}
Let \(\mathcal{G}\) be a lisse \(\overline{\mathbb{Q}}_{\ell}\)-sheaf
on \(\mathbb{A}^{n}_{\overline{\mathbb{F}}_{q}}\), and suppose
\(u\colon \mathbb{A}^{n}_{\overline{\mathbb{F}}_{q}} \to \mathbb{A}_{\overline{\mathbb{F}}_{q}}^{N}\)
is quasi-finite. Then there is a nonempty dense Zariski-open subset
\(\mathscr{U}\) of the dual affine space \(\check{\mathbb{A}}^{N}\)
such that for any hyperplane
\(A \in \mathscr{U}(\overline{\mathbb{F}}_{q})\), the weak Lefschetz
theorem holds for \((u^{-1}A, \mathcal{G})\).
\end{theorem}

\begin{construction}
Consider the closed embedding
\[
\mathbb{A}^n\longrightarrow\mathbb{A}^{n+1},\qquad
(x_1,\ldots,x_n)\longmapsto
(z_{1},\ldots,z_{n},z_{n+1}) = (x_1,\ldots,x_n,x_1x_2+\cdots+x_{n-1}x_n),
\]
and let \(Y \subset \mathbb{A}^n\) be the preimage of a hyperplane
\(a_1z_1 + \cdots + a_nz_n + a_{n+1}z_{n+1} + c = 0\) in
\(\mathbb{A}^{n+1}\). For a general hyperplane we have
\(a_{n+1} \neq 0\), so without loss of generality we may normalize
\(a_{n+1} = 1\); then \(Y = \{P = 0\}\), where
\begin{equation*}
P(x_1,\ldots,x_n)=x_1x_2+\cdots+x_{n-1}x_n+a_1x_1+\cdots+a_nx_n+c.
\end{equation*}
By the weak Lefschetz theorem,
Theorem~\ref{theorem:deligne-weak-lefschetz}, after passing to a
finite extension of \(\mathbb{F}_{q}\), we may choose
\(a_1,\ldots,a_n,c \in \mathbb{F}_q\) so that \(Y\) is smooth and the
weak Lefschetz theorem holds for \((Y, \mathcal{E}^{\boxtimes n})\).
\end{construction}

The key feature of this construction is the quadratic term
\begin{equation*}
x_1x_2 + x_2x_3 + \cdots + x_{n-1}x_n
\end{equation*}
in \(P\).  It is this particular quadratic term that yields the needed
power-saving estimate with asymptotically full square-root
cancellation, as stated in the following elementary lemma.

\begin{lemma}
[\href{https://hexagonmath.org/2611.00789v1}{The Hexagon note}, Lemma 3]
\label{lemma:key}
Let \(n \geq 2\), and let \(\psi \colon \mathbb{F}_q \to \mathbb{C}^{\times}\)
be a nontrivial additive character. For any functions
\(t_1,\ldots,t_n \colon \mathbb{F}_q \to \mathbb{C}\), set
\begin{equation*}
S(t_1,\ldots,t_n) = \sum_{x \in \mathbb{F}_q^n} t_1(x_1) \cdots t_n(x_n)
\cdot \psi\bigl(x_1x_2 + \cdots + x_{n-1}x_n\bigr).
\end{equation*}
Then
\begin{equation*}
|S(t_1 ,\ldots,t_n)| \leq q^{\frac{n+1}{2}} \prod_{i=1}^{n} \|t_i\|_{\infty},
\end{equation*}
where \(\|t\|_{\infty} = \max_{x \in \mathbb{F}_q} |t(x)|\).
\end{lemma}

\begin{remark}
For \(n=2\), Lemma~\ref{lemma:key} reduces to the following
finite-field version of Vinogradov's bilinear bound
\begin{equation*}
|S(t_1,t_2)| \leq q^{\frac{3}{2}} \cdot \|t_1\|_{\infty} \cdot \|t_2\|_{\infty},
\end{equation*}
which saves \(q^{\frac{1}{2}}\) over the trivial bound
\(q^2 \cdot \|t_1\|_{\infty}\cdot \|t_2\|_{\infty}\) for arbitrary bounded functions
\(t_1,t_2\).  See \cite[VI, Problem~14(a)]{vinogradov_elements-number-theory}.

For general \(n\), we have the following  \emph{trivial bound}
\begin{equation*}
|S(t_1,\ldots,t_n)| \leq q^{n}\prod_{i=1}^{n}\|t_i\|_{\infty}
\end{equation*}
for \(S(t_1,\ldots,t_n)\). Thus, for large \(n\),
Lemma~\ref{lemma:key} achieves asymptotic square-root cancellation.
The particular quadratic form
\(x_1x_2 + x_2x_3 + \cdots + x_{n-1}x_n\) is what makes this work.  If
\(n=2k\) and we were instead to use \(k\) disjoint pairs
\(x_1x_2 + x_3x_4 + \cdots + x_{2k-1}x_{2k}\) as the quadratic term,
we would save only \(q^{\frac{k}{2}} = q^{\frac{n}{4}}\), which is
insufficient for the tensor power trick (as we need a saving of
\(q^{\frac{n}{2}-o(n)}\)).
\end{remark}

\begin{proof}[Proof of Proposition~\ref{proposition:main} assuming Lemma~\ref{lemma:key}]
By construction, \(Y = \{P=0\}\) is smooth and satisfies the weak Lefschetz theorem.
Thus, by Lemma~\ref{lem:cayley-trace}, we have
\begin{equation*}
\iota\operatorname{Tr}\bigl(F^{m}, \mathrm{H}_c^n(\mathbb{A}^n;Rj_*j^*(\mathcal{E}^{\boxtimes n}))\bigr)
= (-1)^{n+1}\sum_{h\in\mathbb{F}_{q^m}^{\times}}
\sum_{x\in\mathbb{F}_{q^m}^n}
t(x_1)\cdots t(x_n)\psi_m(hP(x)).
\end{equation*}
Write \(r=\operatorname{rk}\mathcal{E}\). For fixed
\(h\in\mathbb{F}_{q^m}^{\times}\), set \(t_i(x)=t(x)\psi_m(ha_ix)\).
Then \(\|t_i\|_{\infty}\leq r q^{\frac{mw}{2}}\), and the inner sum on the
right-hand side of the displayed equation is
\[
\psi_m(hc)\sum_{x\in\mathbb{F}_{q^m}^n}
t_1(x_1)\cdots t_n(x_n)
\psi_m\bigl(h(x_1x_2+\cdots+x_{n-1}x_n)\bigr).
\]
Applying Lemma~\ref{lemma:key} over \(\mathbb{F}_{q^m}\) with the
nontrivial additive character \(u\mapsto\psi_m(hu)\) and using
\(|\psi_m(hc)|=1\), we see that the inner sum has absolute value at most
\((q^m)^{\frac{n+1}{2}} (r q^{\frac{mw}{2}})^n = r^n q^{\frac{m(nw+n+1)}{2}}\). Summing over \(h\in\mathbb{F}_{q^m}^{\times}\)
gives
\[
\bigl|\iota\operatorname{Tr}\bigl(F^{m},
\mathrm{H}_c^n(\mathbb{A}^n;Rj_*j^*(\mathcal{E}^{\boxtimes n}))\bigr)\bigr|
\leq r^n q^{\frac{m(n(w+1)+3)}{2}}
\]
for all \(m\geq 1\). The proposition therefore follows from
Lemma~\ref{lem:trace-bound}.
\end{proof}

\begin{proof}[Proof of Lemma~\ref{lemma:key}]
We consider two operators on the space of complex-valued functions
on \(\mathbb{F}_q\), equipped with the counting inner product
\(\langle f,g\rangle=\sum_{x\in\mathbb{F}_q} f(x)\overline{g(x)}\)
and its associated \(L^2\)-norm \(\|f\|_2=\langle f,f\rangle^{\frac{1}{2}}\).
The first operator is the (unnormalized) \emph{Fourier transform}
\begin{equation*}
\mathrm{FT}_{\psi}\colon f \mapsto \mathrm{FT}_{\psi}(f),
\qquad
\mathrm{FT}_{\psi}(f)(y) = \sum_{x \in \mathbb{F}_q} f(x) \psi(xy).
\end{equation*}
By Plancherel's identity,
\begin{equation*}
\|\mathrm{FT}_{\psi}(f)\|_{2}=q^{\frac{1}{2}}\|f\|_{2}.
\end{equation*}
The second operator is multiplication by a function \(g\colon \mathbb{F}_q \to \mathbb{C}\):
\begin{equation*}
T_g\colon f \mapsto gf,
\end{equation*}
which satisfies
\begin{equation*}
\|T_g(f)\|_2 = \|gf\|_2 \leq \|g\|_{\infty} \|f\|_2.
\end{equation*}
Using these two operators, we can express the exponential sums
\(S(t_1,\ldots, t_n)\) as iterated inner products:
\begin{align*}
  S(t_1, t_2) &= \langle \mathrm{FT}_{\psi}(t_1), \overline{t}_2 \rangle, \\
  S(t_1,t_2, t_3) &= \langle (\mathrm{FT}_{\psi}\circ T_{t_2} \circ \mathrm{FT}_{\psi})(t_1) ,\overline{t}_3\rangle, \\
  S(t_1,t_2,t_3,t_4) &= \langle (\mathrm{FT}_{\psi} \circ T_{t_3} \circ \mathrm{FT}_{\psi}\circ T_{t_2} \circ \mathrm{FT}_{\psi})(t_1), \overline{t}_4 \rangle,
\end{align*}
and so on. In general, defining \(f_1,\ldots,f_{n-1}\) inductively by
\begin{align*}
  f_1 &= \mathrm{FT}_{\psi}(t_1), \\
  f_{i+1} &= (\mathrm{FT}_{\psi} \circ T_{t_{i+1}})(f_i)
  \qquad (1\leq i\leq n-2),
\end{align*}
we have
\begin{equation*}
\begin{cases}
  S(t_1,\ldots,t_n) = \langle f_{n-1}, \overline{t}_n \rangle, \\[4pt]
  \|f_i\|_2 \leq q^{\frac{i}{2}} \|t_1\|_2 \prod_{j=2}^{i} \|t_j\|_{\infty}
  \qquad (1 \leq i \leq n-1),
\end{cases}
\end{equation*}
where the second bound follows inductively from Plancherel's identity
and the operator norm bound for \(T_{t_{i+1}}\). Applying
Cauchy--Schwarz therefore yields
\begin{align*}
  |S(t_1,\ldots,t_n)|
  &\leq \|f_{n-1}\|_2 \cdot \|t_n\|_2 \\
  &\leq q^{\frac{n-1}{2}}\cdot \|t_1\|_2\cdot \|t_n\|_2 \cdot
       \prod_{i=2}^{n-1}\|t_i\|_{\infty} \\
  &\leq q^{\frac{n+1}{2}}\prod_{i=1}^{n}\|t_i\|_{\infty},
\end{align*}
where the last inequality uses
\(\|t_i\|_2\leq q^{\frac{1}{2}}\|t_i\|_{\infty}\). This completes the
proof.
\end{proof}

\printbibliography
\end{document}